# QUARTER-FRACTION FACTORIAL DESIGNS CONSTRUCTED VIA QUATERNARY CODES[1]

By Frederick K. H. Phoa and Hongquan Xu

*University of California, Los Angeles*

The research of developing a general methodology for the construction of good nonregular designs has been very active in the last decade. Recent research by Xu and Wong [*Statist. Sinica* **17** (2007) 1191–1213] suggested a new class of nonregular designs constructed from quaternary codes. This paper explores the properties and uses of quaternary codes toward the construction of quarter-fraction nonregular designs. Some theoretical results are obtained regarding the aliasing structure of such designs. Optimal designs are constructed under the maximum resolution, minimum aberration and maximum projectivity criteria. These designs often have larger generalized resolution and larger projectivity than regular designs of the same size. It is further shown that some of these designs have generalized minimum aberration and maximum projectivity among all possible designs.

**1. Introduction.** In many scientific researches and investigations, the interests lie in the study of effects of many factors simultaneously. Fractional factorial designs, especially two-level fractional factorial designs, are the most commonly used experimental plans for this type of investigations. Designs that can be constructed through defining relations among factors are called *regular designs*. Any two factorial effects in a regular design are either mutually orthogonal or fully aliased with each other. All other designs that do not possess this kind of defining relationship are called *nonregular designs*.

Regular designs are commonly chosen by the *maximum resolution* criterion [1] and its refinement, the *minimum aberration* criterion [13]. The reader is referred to the books by Wu and Hamada [25] and Mukerjee and Wu [18] for rich results and extensive references.

Received May 2008; revised September 2008.
[1]Supported in part by NSF Grants DMS-05-05728 and DMS-08-06137.
*AMS 2000 subject classification.* 62K15.
*Key words and phrases.* Aliasing index, fractional factorial design, generalized minimum aberration, generalized resolution, nonregular design, projectivity.

The concepts of resolution and aberration have recently been extended to nonregular designs (see [10, 23] and [31]). Tang and Deng [23] showed that generalized minimum aberration designs tend to minimize the contamination of nonnegligible two-factor and higher-order interactions on the estimation of the main effects. Tang [21] provided a projection justification of the generalized minimum aberration criterion, and Cheng, Deng and Tang [7] showed that the generalized minimum aberration criterion is connected with some traditional model-dependent efficiency criteria. For extensions to multi-level nonregular designs, see [9, 17, 26] and [30].

An important and challenging issue is the construction of good nonregular designs. Two simple reasons are: (i) nonregular designs do not have a unified mathematical description, and (ii) there are many more nonregular designs than regular designs. Deng and Tang [11] constructed small generalized minimum aberration designs from Hadamard matrices of order 16, 20 and 24. Tang and Deng [24] constructed generalized resolution designs for 3, 4 and 5 factors and any run size. Li, Deng and Tang [16] searched generalized minimum aberration designs with 20, 24, 28, 32 and 36 runs and up to 6 factors. Xu and Deng [28] searched moment aberration projection designs with 16, 20 and 27 runs. Sun, Li and Ye [20] proposed a sequential algorithm and completely enumerated all 16 and 20-run orthogonal arrays of strength 2. Fang, Zhang and Li [12] proposed an optimization algorithm for construction of generalized minimum aberration designs. Bulutoglu and Margot [3] completely classified some orthogonal arrays of strength 3 up to 56 runs and strength 4 up to 144 runs. All of these algorithmic constructions are limited to small run sizes ($\leq 32$) or small number of factors due to the existence of a large number of designs.

Butler [4] and [5] developed some theoretical results and showed that some existing designs have generalized minimum aberration among all possible designs. Xu [27] constructed several nonregular designs with 32, 64, 128 and 256 runs and 7–16 factors from the Nordstrom and Robinson code, a well-known nonlinear code in coding theory. Tang [22] studied the existence and construction of orthogonal arrays that are robust to nonnegligible two-factor interactions. Stufken and Tang [19] completely classified all two-level orthogonal arrays with $t + 2$ constraints and strength $t$.

In this paper, we consider the construction of two-level nonregular designs via quaternary codes. A quaternary code is a linear space over $Z_4 = \{0, 1, 2, 3\}$, which is the ring of integers modulo 4. Quaternary codes have been successfully used to construct good binary codes in coding theory (see [14]). Xu and Wong [29] first used quaternary codes to construct two-level nonregular designs. They described a systematic procedure for constructing $2^{2n} \times (2^{2n} - 2^n)$ designs and $2^{2n+1} \times (2^{2n+1} - 2^{n+1})$ designs with resolution 3.5 for any $n$, whereas regular designs of the same size have maximum resolution 3 only. They also presented a collection of nonregular designs with 16, 32, 64, 128 and 256 runs and up to 64 factors. Two obvious advantages of using quaternary codes to construct nonregular designs are relatively straightforward construction procedure and simple design presentation. Since the designs are constructed via linear codes over $Z_4$, one can use column indexes to describe these designs. More importantly, many nonregular designs constructed via quaternary codes have better statistical properties than regular designs of the same size in terms of resolution, aberration and projectivity.

The linear structure of a quaternary code makes it possible to analytically study the properties of nonregular designs derived from it. In Section 2, we study the properties of quarter-fraction designs, which can be defined by a generator matrix that consists of an identity matrix and an additional column. It turns out that resolution, wordlength and projectivity can be calculated in terms of the frequency that the numbers 1, 2 and 3 appear in the additional column. Applying these results in Section 3, we construct optimal quarter-fraction designs via quaternary codes under the maximum resolution, minimum aberration and maximum projectivity criteria. These designs are often better than regular designs of the same size in terms of the corresponding criterion. It is well known that a regular minimum aberration design has maximum resolution and maximum projectivity among all regular designs. However, different criteria can lead to different nonregular designs. It turns out that we can often, but not always, find a minimum aberration design that has maximum resolution among all possible quaternary code designs. A minimum aberration design has the same aberration as, and often larger resolution and projectivity than, a regular minimum aberration design. A maximum projectivity design, which often differs from a minimum aberration or maximum resolution design, can have much larger projectivity than a regular minimum aberration design. It is further shown that some of these designs have generalized minimum aberration and maximum projectivity among all possible designs. We present all proofs in Section 4.

The rest of this section introduces notation and definitions. A two-level design $D$, of $N$ runs and $m$ factors, is represented by an $N \times m$ matrix where each row corresponds to a run and each column to a factor, which takes on only two symbols, say $-1$ and $+1$. For $s = \{c_1, c_2, \ldots, c_k\}$, a subset of $k$ columns of $D$, define

$$j_k(s; D) = \sum_{i=1}^{N} c_{i1} \times \cdots \times c_{ik}, \tag{1}$$

where $c_{ij}$ is the $i$th entry of column $c_j$. The $j_k(s; D)$ values are called the *J-characteristics* of design $D$ [10, 21]. It is evident that $|j_k(s; D)| \leq N$.

Following Cheng, Li and Ye [8], we define the *aliasing index* as $\rho_k(s) = \rho_k(s; D) = |j_k(s; D)|/N$, which measures the amount of aliasing among the columns in $s$. It is obvious that $0 \leq \rho_k(s) \leq 1$. When $\rho_k(s) = 1$, the columns in $s$ are fully aliased with each other and form a *complete word* of length $k$. When $0 < \rho_k(s) < 1$, the columns in $s$ are partially aliased with each other and form a *partial word* of length $k$ with aliasing index $\rho_k(s)$. A partial word with aliasing index 1 is a complete word. When $\rho_k(s) = 0$, the columns in $s$ do not form a word.

Suppose that $r$ is the smallest integer such that $\max_{|s|=r} \rho_r(s; D) > 0$, where the maximization is over all subsets of $r$ columns of $D$. The *generalized resolution* [10] of $D$ is defined as

$$(2) \qquad R(D) = r + 1 - \max_{|s|=r} \rho_r(s; D).$$

For $k = 1, \ldots, m$, define

$$(3) \qquad A_k(D) = \sum_{|s|=k} (\rho_k(s; D))^2.$$

The vector $(A_1(D), A_2(D), \ldots, A_m(D))$ is called the generalized wordlength pattern. The *generalized minimum aberration* criterion [30], also called minimum $G_2$-aberration [23], sequentially minimizes the components in the generalized wordlength pattern $A_1(D), A_2(D), \ldots, A_m(D)$. This means that, if two designs have $A_k(D)$ as the first nonequal component in the generalized wordlength pattern, then a design with smaller $A_k(D)$ is preferred.

When restricted to regular designs, generalized resolution, generalized wordlength pattern and generalized minimum aberration reduce to the traditional resolution, wordlength pattern and minimum aberration, respectively. For simplicity, we use resolution, wordlength pattern and minimum aberration for both regular and nonregular designs.

A two-level design $D$ is said to have *projectivity* $p$ [2] if every $p$-factor projection contains a complete $2^p$ factorial design, possibly with some points replicated. It is evident that a regular design of resolution $R = r$ has projectivity $p = r - 1$. Deng and Tang [10] showed that a design with resolution $R > r$ has projectivity $p \geq r$.

## 2. Properties of quarter-fraction designs via quaternary codes.

2.1. *Quaternary codes and binary images.* A quaternary code takes on values from $Z_4 = \{0, 1, 2, 3\}$ (mod. 4). Let $G$ be an $n \times k$ generator matrix over $Z_4$. All possible linear combinations of the rows in $G$ over $Z_4$ form a quaternary linear code, denoted by $C$. The so called *Gray map*, which replaces each element in $Z_4$ with a pair of two symbols, transforms $C$ into a binary code $D = \phi(C)$, which is called the binary image of $C$. For convenience, we use 1 and $-1$ for the two symbols, instead of the 0 and 1 convention for binary codes. Then the Gray map is defined as

$$\phi: 0 \to (1, 1), \qquad 1 \to (1, -1), \qquad 2 \to (-1, -1), \qquad 3 \to (-1, 1).$$

Note that $C$ is a $2^{2n} \times k$ matrix over $Z_4$ and $D$ is a binary $2^{2n} \times 2k$ matrix or a two-level design with $2^{2n}$ runs and $2k$ factors.

2.2. *Designs with $2^{2n}$ runs.* To construct quarter-fraction designs, consider an $n \times (n+1)$ generator matrix $G = (v, I_n)$, where $v$ is an $n \times 1$ column vector over $Z_4$ and $I_n$ is an $n \times n$ identity matrix. Let $D$ be the $2^{2n} \times (2n+2)$ two-level design generated by $G$. It is easy to verify that the identity matrix $I_n$ generates a full $2^{2n} \times 2n$ design; therefore, the properties of $D$ depend on the column $v$ only. Throughout the paper, for $i = 0, 1, 2, 3$, let $f_i$ be the number of times that the number $i$ appears in column $v$. Theorem 1 characterizes, in terms of the frequency $f_i$, the number of words of $D$, their lengths and aliasing indexes.

THEOREM 1. *Consider an $n \times (n+1)$ generator matrix $G = (v, I_n)$. Define $k_1 = f_1 + 2f_2 + f_3 + 1$, $k_2 = 2f_1 + 2f_3 + 2$ and $\rho = 2^{-\lfloor (f_1+f_3)/2 \rfloor}$, where $\lfloor x \rfloor$ is the integer value of $x$. Then, the two-level $2^{2n} \times (2n+2)$ design $D$ generated by $G$ has 1 complete word of length $k_2$ and $2/\rho^2$ partial words of length $k_1$ with aliasing index $\rho$.*

EXAMPLE 1. Consider a generator matrix

$$G = (v \quad I_3) = \begin{pmatrix} 1 & 1 & 0 & 0 \\ 1 & 0 & 1 & 0 \\ 2 & 0 & 0 & 1 \end{pmatrix}.$$

All linear combinations of the three rows of $G$ form a $64 \times 4$ linear code $C$ over $Z_4$. Applying the Gray map, a $64 \times 8$ binary image $D = \phi(C)$ is obtained; see Table 1. According to Theorem 1, design $D$ has 1 complete word of length $k_2 = 6$ and 8 partial words of length $k_1 = 5$ with aliasing index $\rho = 0.5$. It is easy to verify that the first six columns form a complete word and that columns $(a, b, c, 7, 8)$ form a partial word with aliasing index 0.5, where $a = 1$ or 2, $b = 3$ or 4 and $c = 5$ or 6. Therefore, by definitions (2) and (3), the resolution of $D$ is 5.5, and the wordlength pattern of $D$ is $A_5(D) = 2$, $A_6(D) = 1$ and $A_i(D) = 0$ for $i \neq 5, 6$.

For ease of presentation, we say that the $i$th identity column of $I_n$ in $G = (v, I_n)$ is "associated with" number $z$ if the $i$th element of $v$ is $z$, where $z = 0, 1, 2$ or 3. We also refer to a column of $D$ as associated with number $z$ if it is one of the two columns generated by an identity column that is associated with number $z$. Further, we refer to the two columns generated by $v$ as associated with vector $v$. For example, the first two columns of $D$ in Table 1 are associated with vector $v$, columns 3 to 6 are associated with number 1, and the last two columns are associated with number 2.

Now, we can describe more precisely about the words of $D$ in Theorem 1. The complete word of $D$ consists of all columns associated with vector $v$ and numbers 1 and 3. Each partial word consists of all columns associated with number 2, one of the columns associated with vector $v$ and each number 1 and 3. Furthermore, the columns associated with number 0 do not appear in any word.

Recall that a regular design has only complete words. Corollary 1 provides a sufficient and necessary condition for $D$ to be a regular design.

COROLLARY 1. *Design $D$ is regular if and only if $f_1 + f_3 \leq 1$.*

It is straightforward to complete the resolution of $D$ according to the definition (2) and Theorem 1.

TABLE 1
*A quaternary code $C$ and its binary image $D$*

| | Code $C$ | | | | | Design $D$ | | | | | | | |
|---|---|---|---|---|---|---|---|---|---|---|---|---|---|
| Run | 1 | 2 | 3 | 4 | Run | 1 | 2 | 3 | 4 | 5 | 6 | 7 | 8 |
| 1 | 0 | 0 | 0 | 0 | 1 | 1 | 1 | 1 | 1 | 1 | 1 | 1 | 1 |
| 2 | 1 | 1 | 0 | 0 | 2 | 1 | −1 | 1 | −1 | 1 | 1 | 1 | 1 |
| 3 | 2 | 2 | 0 | 0 | 3 | −1 | −1 | −1 | −1 | 1 | 1 | 1 | 1 |
| 4 | 3 | 3 | 0 | 0 | 4 | −1 | 1 | −1 | 1 | 1 | 1 | 1 | 1 |
| 5 | 1 | 0 | 1 | 0 | 5 | 1 | −1 | 1 | 1 | 1 | −1 | 1 | 1 |
| 6 | 2 | 1 | 1 | 0 | 6 | −1 | −1 | 1 | −1 | 1 | −1 | 1 | 1 |
| 7 | 3 | 2 | 1 | 0 | 7 | −1 | 1 | −1 | −1 | 1 | −1 | 1 | 1 |
| 8 | 0 | 3 | 1 | 0 | 8 | 1 | 1 | −1 | 1 | 1 | −1 | 1 | 1 |
| 9 | 2 | 0 | 2 | 0 | 9 | −1 | −1 | 1 | 1 | −1 | −1 | 1 | 1 |
| 10 | 3 | 1 | 2 | 0 | 10 | −1 | 1 | 1 | −1 | −1 | −1 | 1 | 1 |
| 11 | 0 | 2 | 2 | 0 | 11 | 1 | 1 | −1 | −1 | −1 | −1 | 1 | 1 |
| 12 | 1 | 3 | 2 | 0 | 12 | 1 | −1 | −1 | 1 | −1 | −1 | 1 | 1 |
| 13 | 3 | 0 | 3 | 0 | 13 | −1 | 1 | 1 | 1 | −1 | 1 | 1 | 1 |
| 14 | 0 | 1 | 3 | 0 | 14 | 1 | 1 | 1 | −1 | −1 | 1 | 1 | 1 |
| 15 | 1 | 2 | 3 | 0 | 15 | 1 | −1 | −1 | −1 | −1 | 1 | 1 | 1 |
| 16 | 2 | 3 | 3 | 0 | 16 | −1 | −1 | −1 | 1 | −1 | 1 | 1 | 1 |
| 17 | 2 | 0 | 0 | 1 | 17 | −1 | −1 | 1 | 1 | 1 | 1 | 1 | −1 |
| 18 | 3 | 1 | 0 | 1 | 18 | −1 | 1 | 1 | −1 | 1 | 1 | 1 | −1 |
| 19 | 0 | 2 | 0 | 1 | 19 | 1 | 1 | −1 | −1 | 1 | 1 | 1 | −1 |
| 20 | 1 | 3 | 0 | 1 | 20 | 1 | −1 | −1 | 1 | 1 | 1 | 1 | −1 |
| 21 | 3 | 0 | 1 | 1 | 21 | −1 | 1 | 1 | 1 | 1 | −1 | 1 | −1 |
| 22 | 0 | 1 | 1 | 1 | 22 | 1 | 1 | 1 | −1 | 1 | −1 | 1 | −1 |
| 23 | 1 | 2 | 1 | 1 | 23 | 1 | −1 | −1 | −1 | 1 | −1 | 1 | −1 |
| 24 | 2 | 3 | 1 | 1 | 24 | −1 | −1 | −1 | 1 | 1 | −1 | 1 | −1 |
| 25 | 0 | 0 | 2 | 1 | 25 | 1 | 1 | 1 | 1 | −1 | −1 | 1 | −1 |

TABLE 1
*(Continued)*

| | Code $C$ | | | | | Design $D$ | | | | | | | |
|---|---|---|---|---|---|---|---|---|---|---|---|---|---|
| Run | 1 | 2 | 3 | 4 | Run | 1 | 2 | 3 | 4 | 5 | 6 | 7 | 8 |
| 26 | 1 | 1 | 2 | 1 | 26 | 1 | $-1$ | 1 | $-1$ | $-1$ | $-1$ | 1 | $-1$ |
| 27 | 2 | 2 | 2 | 1 | 27 | $-1$ | $-1$ | $-1$ | $-1$ | $-1$ | $-1$ | 1 | $-1$ |
| 28 | 3 | 3 | 2 | 1 | 28 | $-1$ | 1 | $-1$ | 1 | $-1$ | $-1$ | 1 | $-1$ |
| 29 | 1 | 0 | 3 | 1 | 29 | 1 | $-1$ | 1 | 1 | $-1$ | 1 | 1 | $-1$ |
| 30 | 2 | 1 | 3 | 1 | 30 | $-1$ | $-1$ | 1 | $-1$ | $-1$ | 1 | 1 | $-1$ |
| 31 | 3 | 2 | 3 | 1 | 31 | $-1$ | 1 | $-1$ | $-1$ | $-1$ | 1 | 1 | $-1$ |
| 32 | 0 | 3 | 3 | 1 | 32 | 1 | 1 | $-1$ | 1 | $-1$ | 1 | 1 | $-1$ |
| 33 | 0 | 0 | 0 | 2 | 33 | 1 | 1 | 1 | 1 | 1 | 1 | $-1$ | $-1$ |
| 34 | 1 | 1 | 0 | 2 | 34 | 1 | $-1$ | 1 | $-1$ | 1 | 1 | $-1$ | $-1$ |
| 35 | 2 | 2 | 0 | 2 | 35 | $-1$ | $-1$ | $-1$ | $-1$ | 1 | 1 | $-1$ | $-1$ |
| 36 | 3 | 3 | 0 | 2 | 36 | $-1$ | 1 | $-1$ | 1 | 1 | 1 | $-1$ | $-1$ |
| 37 | 1 | 0 | 1 | 2 | 37 | 1 | $-1$ | 1 | 1 | 1 | $-1$ | $-1$ | $-1$ |
| 38 | 2 | 1 | 1 | 2 | 38 | $-1$ | $-1$ | 1 | $-1$ | 1 | $-1$ | $-1$ | $-1$ |
| 39 | 3 | 2 | 1 | 2 | 39 | $-1$ | 1 | $-1$ | $-1$ | 1 | $-1$ | $-1$ | $-1$ |
| 40 | 0 | 3 | 1 | 2 | 40 | 1 | 1 | $-1$ | 1 | 1 | $-1$ | $-1$ | $-1$ |
| 41 | 2 | 0 | 2 | 2 | 41 | $-1$ | $-1$ | 1 | 1 | $-1$ | $-1$ | $-1$ | $-1$ |
| 42 | 3 | 1 | 2 | 2 | 42 | $-1$ | 1 | 1 | $-1$ | $-1$ | $-1$ | $-1$ | $-1$ |
| 43 | 0 | 2 | 2 | 2 | 43 | 1 | 1 | $-1$ | $-1$ | $-1$ | $-1$ | $-1$ | $-1$ |
| 44 | 1 | 3 | 2 | 2 | 44 | 1 | $-1$ | $-1$ | 1 | $-1$ | $-1$ | $-1$ | $-1$ |
| 45 | 3 | 0 | 3 | 2 | 45 | $-1$ | 1 | 1 | 1 | $-1$ | 1 | $-1$ | $-1$ |
| 46 | 0 | 1 | 3 | 2 | 46 | 1 | 1 | 1 | $-1$ | $-1$ | 1 | $-1$ | $-1$ |
| 47 | 1 | 2 | 3 | 2 | 47 | 1 | $-1$ | $-1$ | $-1$ | $-1$ | 1 | $-1$ | $-1$ |
| 48 | 2 | 3 | 3 | 2 | 48 | $-1$ | $-1$ | $-1$ | 1 | $-1$ | 1 | $-1$ | $-1$ |
| 49 | 2 | 0 | 0 | 3 | 49 | $-1$ | $-1$ | 1 | 1 | 1 | 1 | $-1$ | 1 |
| 50 | 3 | 1 | 0 | 3 | 50 | $-1$ | 1 | 1 | $-1$ | 1 | 1 | $-1$ | 1 |
| 51 | 0 | 2 | 0 | 3 | 51 | 1 | 1 | $-1$ | $-1$ | 1 | 1 | $-1$ | 1 |
| 52 | 1 | 3 | 0 | 3 | 52 | 1 | $-1$ | $-1$ | 1 | 1 | 1 | $-1$ | 1 |
| 53 | 3 | 0 | 1 | 3 | 53 | $-1$ | 1 | 1 | 1 | 1 | $-1$ | $-1$ | 1 |
| 54 | 0 | 1 | 1 | 3 | 54 | 1 | 1 | 1 | $-1$ | 1 | $-1$ | $-1$ | 1 |
| 55 | 1 | 2 | 1 | 3 | 55 | 1 | $-1$ | $-1$ | $-1$ | 1 | $-1$ | $-1$ | 1 |
| 56 | 2 | 3 | 1 | 3 | 56 | $-1$ | $-1$ | $-1$ | 1 | 1 | $-1$ | $-1$ | 1 |
| 57 | 0 | 0 | 2 | 3 | 57 | 1 | 1 | 1 | 1 | $-1$ | $-1$ | $-1$ | 1 |
| 58 | 1 | 1 | 2 | 3 | 58 | 1 | $-1$ | 1 | $-1$ | $-1$ | $-1$ | $-1$ | 1 |
| 59 | 2 | 2 | 2 | 3 | 59 | $-1$ | $-1$ | $-1$ | $-1$ | $-1$ | $-1$ | $-1$ | 1 |
| 60 | 3 | 3 | 2 | 3 | 60 | $-1$ | 1 | $-1$ | 1 | $-1$ | $-1$ | $-1$ | 1 |
| 61 | 1 | 0 | 3 | 3 | 61 | 1 | $-1$ | 1 | 1 | $-1$ | 1 | $-1$ | 1 |
| 62 | 2 | 1 | 3 | 3 | 62 | $-1$ | $-1$ | 1 | $-1$ | $-1$ | 1 | $-1$ | 1 |
| 63 | 3 | 2 | 3 | 3 | 63 | $-1$ | 1 | $-1$ | $-1$ | $-1$ | 1 | $-1$ | 1 |
| 64 | 0 | 3 | 3 | 3 | 64 | 1 | 1 | $-1$ | 1 | $-1$ | 1 | $-1$ | 1 |

COROLLARY 2. *The resolution of $D$ is $k_2$ if $k_1 \geq k_2$, or $k_1 + 1 - \rho$ otherwise.*

According to the definition (3), when summing up $2/\rho^2$ partial words of length $k_1$ with aliasing index $\rho$, we get $A_{k_1}(D) = 2$. Corollary 3 specifies the wordlength pattern of $D$.

COROLLARY 3. *The wordlength pattern of $D$ is:*
(a) *If $k_1 \neq k_2$, then $A_{k_1}(D) = 2$, $A_{k_2}(D) = 1$ and $A_i(D) = 0$ for $i \neq k_1, k_2$;*
(b) *If $k_1 = k_2 = k$, then $A_k(D) = 3$ and $A_i(D) = 0$ for $i \neq k$.*

Next, we consider the projectivity of design $D$ generated by $G = (v, I_n)$. Theorem 1 suggests that there is a complete word of length $k_2 = 2(f_1 + f_3) + 2$. This implies that the projectivity of $D$ is, at most, $2(f_1 + f_3) + 1$. The next theorem states that the projectivity of $D$ is indeed $2(f_1 + f_3) + 1$ if $f_2 > 0$.

THEOREM 2. *Suppose that $D$ is the two-level $2^{2n} \times (2n+2)$ design generated by $G = (v, I_n)$:*
(a) *If $f_2 > 0$, the projectivity of $D$ is $2(f_1 + f_3) + 1$;*
(b) *If $f_2 = 0$ and $f_1 + f_3 > 0$, the projectivity of $D$ is $2(f_1 + f_3) - 1$.*

Theorem 2 implies that the projectivity of $D$ is not affected by the partial words. As an example, consider design $D$ in Example 1. Theorem 2 suggests that the projectivity of $D$ is 5. This can be verified directly.

2.3. *Designs with $2^{2n-1}$ runs.* Design $D$, generated by $G = (v, I_n)$, has $2^{2n}$ runs and $2n+2$ factors. To construct quarter-fraction designs with $2^{2n-1}$ runs, we use the half fraction method, which works as follows. Choose any column of $D$ as a branching column, which divides $D$ into two half-fractions according to the symbols of the branching column. Deleting the branching column yields two $2^{2n-1} \times (2n+1)$ designs. It is easy to verify that the two half-fractions of $D$ are equivalent. However, the properties of the half-fractions depend on the branching column, which are characterized in Theorem 3.

THEOREM 3. *Suppose that $D$ is the two-level $2^{2n} \times (2n+2)$ design generated by $G = (v, I_n)$ and that $D'$ is a half-fraction of $D$. Define $k_1$, $k_2$ and $\rho$ as in Theorem 1:*
(a) *If the branching column is associated with number 1 or 3, $D'$ has 1 complete word of length $k_2 - 1$, $1/\rho^2$ partial words of length $k_1$ with aliasing index $\rho$ and $1/\rho^2$ partial words of length $k_1 - 1$ with aliasing index $\rho$;*
(b) *If the branching column is associated with number 2, $D'$ has 1 complete word of length $k_2$ and $2/\rho^2$ partial words of length $k_1 - 1$ with aliasing index $\rho$.*

It is easy to verify that, if the branching column is associated with vector $v$, this is identical to case (a) when $f_1 + f_3 > 0$ or case (b) when $f_1 + f_3 = 0$ and $f_2 > 0$. If the branching column is associated with number 0, $D'$ and $D$ share the same words because the branching column does not appear in any word of $D$.

The following four corollaries summarize the resolution and wordlength pattern of $D'$ for cases (a) and (b), separately.

COROLLARY 4. *The resolution of $D'$ derived in Theorem 3(a) is $k_2 - 1$ if $k_1 \geq k_2$, or $k_1 - \rho$ otherwise.*

COROLLARY 5. *The wordlength pattern of $D'$ derived in Theorem 3(a) is:*

(a) *If $k_1 = k_2 = k$, then $A_{k-1}(D') = 2$, $A_k(D') = 1$ and $A_i(D') = 0$ for $i \neq k-1, k$;*

(b) *If $k_1 = k_2 - 1 = k$, then $A_{k-1}(D') = 1$, $A_k(D') = 2$ and $A_i(D') = 0$ for $i \neq k-1, k$;*

(c) *If $k_1 \neq k_2$ or $k_2 - 1$, then $A_{k_1-1}(D') = A_{k_2-1}(D') = A_{k_1}(D') = 1$ and $A_i(D') = 0$ for $i \neq k_1 - 1, k_1, k_2 - 1$.*

COROLLARY 6. *The resolution of $D'$ derived in Theorem 3(b) is $k_2$ if $k_1 - 1 \geq k_2$, or $k_1 - \rho$ otherwise.*

COROLLARY 7. *The wordlength pattern of $D'$ derived in Theorem 3(b) is:*

(a) *If $k_1 - 1 \neq k_2$, then $A_{k_1-1}(D') = 2$, $A_{k_2}(D') = 1$ and $A_i(D') = 0$ for $i \neq k_1 - 1, k_2$;*

(b) *If $k_1 - 1 = k_2 = k$, then $A_k(D') = 3$ and $A_i(D') = 0$ for $i \neq k$.*

The next theorem summarizes the projectivity of a half-fraction of $D$.

THEOREM 4. *Suppose that $D$ is the two-level $2^{2n} \times (2n+2)$ design generated by $G = (v, I_n)$ and that $D'$ is a half-fraction of $D$:*

(a) *If $f_2 > 0$, $f_1 + f_3 > 0$ and the branching column is associated with number 1 or 3, the projectivity of $D'$ is $2(f_1 + f_3)$;*

(b) *If $f_2 = 0$, $f_1 + f_3 > 0$ and the branching column is associated with number 1 or 3, the projectivity of $D'$ is $2(f_1 + f_3) - 2$;*

(c) *If $f_2 > 1$ and the branching column is associated with number 2, the projectivity of $D'$ is $2(f_1 + f_3) + 1$;*

(d) *If $f_2 = 1$ and the branching column is associated with number 2, the projectivity of $D'$ is $2(f_1 + f_3)$.*

Comparing with Theorem 2, we observe that the projectivity of $D'$ is equal to the projectivity of $D$ for case (c), whereas the projectivity of $D'$ is equal to the projectivity of $D$ minus one for all other cases.

EXAMPLE 2. Consider half-fractions of $D$ in Table 1. If one of the first six columns is chosen as the branching column, we obtain a $32 \times 7$ design $D'$ with resolution 4.5 and wordlength patterns $A_4(D') = 1$, $A_5(D') = 2$ and $A_i(D') = 0$ for $i \neq 4, 5$. Design $D'$ has 1 complete word of length 5, 4 partial words of length 5 with aliasing index 0.5 and 4 partial words of length 4 with aliasing index 0.5. For example, if the first column is chosen as the branching column, then columns 2 to 6 form a complete word and columns $(b, c, 7, 8)$ and $(2, b, c, 7, 8)$ form a partial word with aliasing index 0.5, where $b = 3$ or 4 and $c = 5$ or 6. If one of the last two columns is chosen as the branching column, we obtain a $32 \times 7$ design $D'$ with resolution 4.5 and wordlength patterns $A_4(D') = 2$, $A_6(D') = 1$ and $A_i(D') = 0$ for $i \neq 4, 6$. Design $D'$ has 1 complete word of length 6 and 8 partial words of length 4 with aliasing index 0.5. Finally, according to Theorem 4, any half-fraction of $D$ has projectivity 4, which can be verified directly.

**3. Optimal quarter-fraction designs.** In this section, we apply the theory developed in the previous section to construct optimal designs under the maximum resolution, minimum aberration and maximum projectivity criteria. As shown below, different criteria can lead to different optimal designs.

3.1. *Designs with $2^{2n}$ runs.* Applying Theorem 1, we have the following results regarding maximum resolution and minimum aberration designs.

THEOREM 5. *Among all $2^{2n} \times (2n+2)$ designs generated by $G = (v, I_n)$:*
(a) *If $n = 3k - 1, k \geq 1$, then a design $D$ defined by $f_1 + f_3 = 2k - 1$ and $f_2 = k$ has maximum resolution $4k$;*
(b) *If $n = 3k, k \geq 1$, then a design $D$ defined by $f_1 + f_3 = 2k$ and $f_2 = k$ has maximum resolution $4k + 2 - 2^{-k}$;*
(c) *If $n = 3k + 1, k \geq 1$, then a design $D$ defined by $f_1 + f_3 = 2k + 1$ and $f_2 = k$ has maximum resolution $4k + 3 - 2^{-k}$.*

THEOREM 6. *Among all $2^{2n} \times (2n+2)$ designs generated by $G = (v, I_n)$:*
(a) *If $n = 3k - 1, k \geq 1$, then a design $D$ defined by $f_1 + f_3 = 2k - 1$ and $f_2 = k$ has minimum aberration and its wordlength pattern is $A_{4k}(D) = 3$;*
(b) *If $n = 3k, k \geq 1$, then a design $D$ defined by $f_1 + f_3 = 2k$ and $f_2 = k$ has minimum aberration and its wordlength pattern is $A_{4k+1}(D) = 2$ and $A_{4k+2}(D) = 1$;*
(c) *If $n = 3k + 1, k \geq 1$, then a design $D$ defined by $f_1 + f_3 = 2k$ and $f_2 = k + 1$ has minimum aberration and its wordlength pattern is $A_{4k+2}(D) = 1$ and $A_{4k+3}(D) = 2$.*

When $n = 3k - 1$ or $3k$, the minimum aberration design in Theorem 6 coincides with the maximum resolution design in Theorem 5; however, when $n = 3k + 1$, the minimum aberration design differs from the maximum resolution design.

Applying Theorem 2, we have the following result regarding maximum projectivity designs.

THEOREM 7. *Among all $2^{2n} \times (2n+2)$ designs generated by $G = (v, I_n)$, a design $D$ defined by $f_1 + f_3 = n - 1$ and $f_2 = 1$ has maximum projectivity $2n - 1$, and so does a design $D$ defined by $f_1 + f_3 = n$ and $f_2 = 0$.*

The maximum projectivity designs in Theorem 7 are different from designs in Theorems 5 and 6 when $n > 4$. According to Corollary 2, a design defined by $f_1 + f_3 = n - 1$ and $f_2 = 1$ has resolution $n + 3 - 2^{-\lfloor (n-1)/2 \rfloor}$ for $n \geq 2$, and a design defined by $f_1 + f_3 = n$ and $f_2 = 0$ has resolution $n + 2 - 2^{-\lfloor n/2 \rfloor}$; therefore, the former design is recommended.

3.2. *Designs with $2^{2n-1}$ runs.* To find optimal designs with $2^{2n-1}$ runs, we consider all possible designs generated by $G = (v, I_n)$ and all possible half-fractions. It turns out that it is sufficient to consider only half-fractions of the minimum aberration designs in Theorem 6 and the maximum projectivity designs in Theorem 7.

THEOREM 8. *Suppose that $D'$ is a half-fraction of a design $D$ given in Theorem 6. Among all $2^{2n-1} \times (2n+1)$ designs that are half-fractions of designs generated by $G = (v, I_n)$, $D'$ has maximum resolution and minimum aberration:*

*(a) If $n = 3k - 1, k \geq 1$, and the branching column is associated with number 2. The resolution of $D'$ is $4k - 2^{-(k-1)}$ and the wordlength pattern is $A_{4k-1}(D') = 2$ and $A_{4k}(D') = 1$;*

*(b) If $n = 3k, k \geq 1$, and the branching column is associated with number 1. The resolution of $D'$ is $4k + 1 - 2^{-k}$ and the wordlength pattern is $A_{4k}(D') = 1$ and $A_{4k+1}(D') = 2$;*

*(c) If $n = 3k + 1, k \geq 1$, and the branching column is associated with number 2. The resolution of $D'$ is $4k + 2$ and the wordlength pattern is $A_{4k+2}(D') = 3$.*

THEOREM 9. *Any half-fraction of a design $D$ in Theorem 7 has maximum projectivity $2n - 2$ among all $2^{2n-1} \times (2n+1)$ designs that are half-fractions of designs generated by $G = (v, I_n)$.*

TABLE 2
*Optimal quarter-fraction designs*

| Design | Quaternary-code designs | | | | | Regular | |
| --- | --- | --- | --- | --- | --- | --- | --- |
| | Criterion | $v^T$ | WLP | $R$ | pr | $R$ | pr |
| $2^{6-2}$ | $r, a, p$ | [12] | $A_4 = 3$ | 4.0 | 3 | 4 | 3 |
| $2^{7-2}$ | $r, a, p$ | $[112]f$ | $A_4 = 1, A_5 = 2$ | 4.5 | 4 | 4 | 3 |
| $2^{8-2}$ | $r, a, p$ | [112] | $A_5 = 2, A_6 = 1$ | 5.5 | 5 | 5 | 4 |
| $2^{9-2}$ | $r, a$ | $[1122]l$ | $A_6 = 3$ | 6.0 | 5 | 6 | 5 |
| | $p$ | $[1112]f$ | $A_5 = 1, A_6 = 2$ | 5.5 | 6 | | |
| $2^{10-2}$ | $r, p$ | [1112] | $A_6 = 2, A_8 = 1$ | 6.5 | 7 | 6 | 5 |
| | $a$ | [1122] | $A_6 = 1, A_7 = 2$ | 6.0 | 5 | | |
| $2^{11-2}$ | $r, a$ | $[11122]l$ | $A_7 = 2, A_8 = 1$ | 7.5 | 7 | 7 | 6 |
| | $p$ | $[11112]f$ | $A_6 = A_7 = A_9 = 1$ | 6.75 | 8 | | |
| $2^{12-2}$ | $r, a$ | [11122] | $A_8 = 3$ | 8.0 | 7 | 8 | 7 |
| | $p$ | [11112] | $A_7 = 2, A_{10} = 1$ | 7.75 | 9 | | |
| $2^{13-2}$ | $r, a$ | $[111122]f$ | $A_8 = 1, A_9 = 2$ | 8.75 | 8 | 8 | 7 |
| | $p$ | $[111112]f$ | $A_7 = A_8 = A_{11} = 1$ | 7.75 | 10 | | |
| $2^{14-2}$ | $r, a$ | [111122] | $A_9 = 2, A_{10} = 1$ | 9.75 | 9 | 9 | 8 |
| | $p$ | [111112] | $A_8 = 2, A_{12} = 1$ | 8.75 | 11 | | |
| $2^{15-2}$ | $r, a$ | $[1111222]l$ | $A_{10} = 3$ | 10.0 | 9 | 10 | 9 |
| | $p$ | $[1111112]f$ | $A_8 = A_9 = A_{13} = 1$ | 8.875 | 12 | | |
| $2^{16-2}$ | $r$ | [1111122] | $A_{10} = 2, A_{12} = 1$ | 10.75 | 11 | 10 | 9 |
| | $a$ | [1111222] | $A_{10} = 1, A_{11} = 2$ | 10.0 | 9 | | |
| | $p$ | [1111112] | $A_9 = 2, A_{14} = 1$ | 9.875 | 13 | | |

3.3. *Table of designs.* For easy reference, we provide some optimal designs and their properties in Table 2. Following the convention on regular designs, we use the notation $2^{m-2}$ to represent a quarter-fraction design with $m$ factors and $2^{m-2}$ runs. The second column of Table 2 specifies the three optimality criteria: maximum resolution (r), minimum aberration (a) and maximum projectivity (p). The third column is the vector $v$ in the generator matrix $G = (v, I_n)$ and the letter at the end denotes the branching column, which is either the first ($f$) or last ($l$) column. The first column is associated with vector $v$, while the last column is associated with number 2. Choosing the first column or a column associated with number 1 as the branching column yields an equivalent design. The next three columns, under the category of "quaternary-code designs," are the wordlength pattern (WLP), resolution (R) and projectivity (pr) of the design generated by $G = (v, I_n)$. The last two columns, under the category of "regular," are the resolution and projectivity of a regular minimum aberration design with the same size.

Table 2 shows that the maximum resolution designs and the minimum aberration designs are similar, but they often differ from the maximum projectivity designs. Specifically, the "r" design coincides with the "a" design when $m \neq 6k + 4$, $k > 0$, whereas the "$p$" design differs from the "$r$" or "$a$" design when $m = 9$ or $m > 10$.

According to Corollary 1, all designs in Table 2 are nonregular designs, except for design $2^{6-2}$, which is equivalent to the regular minimum aberration design. Design $2^{8-2}$ is considered in Example 1 and given explicitly in Table 1. Design $2^{7-2}$ is a half-fraction of design $2^{8-2}$ and illustrated in Example 2.

It is of great interest to compare the quaternary-code designs with regular minimum aberration $2^{m-2}$ designs, which were given by Chen and Wu [6]. A regular minimum aberration $2^{m-2}$ design has resolution $R = \lfloor 2m/3 \rfloor$, projectivity $R - 1$ and wordlength pattern $A_R = 3R - 2m + 3$ and $A_{R+1} = 2m - 3R$. All of the "$r$" designs in Table 2 have the same or larger resolution as regular minimum aberration designs; in particular, when $m = 3k + 1$ or $3k + 2$, all of the "$r$" designs have larger resolution and, therefore, larger projectivity. All of the "$a$" designs have the same wordlength pattern as regular minimum aberration designs and have the same or larger resolution and projectivity. Indeed, Xu [27] showed that regular minimum aberration $2^{m-2}$ designs have minimum aberration among all possible designs. Except for design $2^{6-2}$, all of the "$p$" designs have higher projectivity than regular minimum aberration designs, but they may have smaller resolution. Indeed, all of the "$p$" designs have maximum projectivity among all possible designs. The next theorem summarizes these results.

THEOREM 10. (a) *The designs given in Theorems 6 and 8 have minimum aberration among all possible designs.*

(b) *The designs given in Theorems 7 and 9 have maximum projectivity among all possible designs.*

It is of interest to know whether the designs given in Theorems 5 and 8 have maximum resolution among all possible designs. We do not have an answer yet. The compete catalogs of [3, 20] suggest that designs $2^{6-2}$, $2^{7-2}$ and $2^{8-2}$ given in Table 2 have maximum resolution among all possible designs. This can also be verified analytically using Proposition 2 of Deng and Tang [10].

Another interesting question is whether the optimality results can be extended to 1/16 fraction designs by using a generator matrix which consists of an identity matrix plus two columns. This is much more complicated due to the fact that we have to deal with 16 level combinations of the two extra columns. We are investigating this problem.

**4. Proofs.** Some lemmas are introduced in order to prove the theorems.

4.1. *Some lemmas.* Consider an $n \times (n+1)$ generator matrix $G_n = (v_n, I_n)$, where $v_n$ is an $n \times 1$ column vector over $Z_4$ and $I_n$ is an $n \times n$ identity matrix. Let $D_n$ be the $2^{2n} \times (2n+2)$ binary design generated by $G_n$.

Let $v_{n-1}$ be the vector consisting of the first $n-1$ components of $v_n$, and let $D_{n-1}$ be the $2^{2n-2} \times 2n$ binary design generated by the $(n-1) \times n$ generator matrix $G_{n-1} = (v_{n-1}, I_{n-1})$. Denote $D_{n-1} = (a, b, E)$, where $a$ and $b$ are column vectors generated by $v_{n-1}$, and $E$ is a $2^{(2n-2)} \times (2n-2)$ full factorial generated by $I_{n-1}$.

We can express $D_n$ in terms of $D_{n-1}$, depending on the last component of $v_n$, which is denoted by $z$. It is trivial for $z = 0$. It is obvious that $z = 1$ and $z = 3$ produce an equivalent design. Therefore, it is sufficient to consider only $z = 1$ or 2.

When $z = 1$, $D_n$ can be expressed as follows, up to row permutations

$$(4) \qquad D_n = \begin{pmatrix} a & b & E & 1 & 1 \\ b & -a & E & 1 & -1 \\ -a & -b & E & -1 & -1 \\ -b & a & E & -1 & 1 \end{pmatrix},$$

where **1** is a vector of ones. From this expression and the definition (1), we establish the connection between the $J$-characteristics of $D_n$ and $D_{n-1}$. Note that the column indexes of $D_n$ are $\{1, 2, \ldots, 2n+2\}$ and of $D_{n-1}$ are $\{1, 2, \ldots, 2n\}$. For clarification, the $s$ in the notation $j_k(s; D)$ refers to a subset of column indexes of $D$, and we omit $k$ when it is not important.

LEMMA 1. *Suppose that the last component of $v_n$ is 1. For any subset $e \subset \{3, 4, \ldots, 2n\}$:*
(a) $j(\{1, 2n+1\} \cup e; D_n) = j(\{2, 2n+2\} \cup e; D_n) = 2j(\{1\} \cup e; D_{n-1}) + 2j(\{2\} \cup e; D_{n-1})$;
(b) $j(\{1, 2n+2\} \cup e; D_n) = -j(\{2, 2n+1\} \cup e; D_n) = 2j(\{1\} \cup e; D_{n-1}) - 2j(\{2\} \cup e; D_{n-1})$;
(c) $j(\{1, 2, 2n+1, 2n+2\} \cup e; D_n) = 4j(\{1, 2\} \cup e; D_{n-1})$;
(d) $j(s \cup e; D_n) = 0$ for $s = \{1\}, \{2\}, \{2n+1\}, \{2n+2\}, \{1, 2\}, \{2n+1, 2n+2\}, \{1, 2, 2n+1\}, \{1, 2, 2n+2\}, \{1, 2n+1, 2n+2\}$, or $\{2, 2n+1, 2n+2\}$.

When $z = 2$, $D_n$ can be expressed as follows, up to row permutations,

$$(5) \qquad D_n = \begin{pmatrix} a & b & E & 1 & 1 \\ -a & -b & E & 1 & -1 \\ a & b & E & -1 & -1 \\ -a & -b & E & -1 & 1 \end{pmatrix}.$$

From this expression and the definition (1), we establish the connection between the $J$-characteristics of $D_n$ and $D_{n-1}$.

LEMMA 2. *Suppose that the last component of $v_n$ is 2. For any subset $e \subset \{3, 4, \ldots, 2n\}$:*
  (a) $j(\{1, 2n+1, 2n+2\} \cup e; D_n) = 4j(\{1\} \cup e; D_{n-1})$;
  (b) $j(\{2, 2n+1, 2n+2\} \cup e; D_n) = 4j(\{2\} \cup e; D_{n-1})$;
  (c) $j(\{1, 2\} \cup e; D_n) = 4j(\{1, 2\} \cup e; D_{n-1})$;
  (d) $j(s \cup e; D_n) = 0$ for $s = \{1\}, \{2\}, \{2n+1\}, \{2n+2\}, \{1, 2n+1\}, \{1, 2n+2\}, \{2, 2n+1\}, \{2, 2n+2\}, \{2n+1, 2n+2\}, \{1, 2, 2n+1\}, \{1, 2, 2n+2\}$, or $\{1, 2, 2n+1, 2n+2\}$.

The next result describes the partial words of $D_n$ and their $J$-characteristics.

LEMMA 3. *Suppose that $v_n$ is a vector of $n$ 1's. For $l = 1, 2$, let $s_l = \{l, x_2, \ldots, x_{n+1}\}$ where $x_i = 2i - 1$ or $2i$ for $i = 2, \ldots, n+1$:*
  (a) *If $n = 2t + 1$, either $j_{n+1}(s_1; D_n) = 0$ and $|j_{n+1}(s_2; D_n)| = 2^{3t+2}$ or $|j_{n+1}(s_1; D_n)| = 2^{3t+2}$ and $j_{n+1}(s_2; D_n) = 0$;*
  (b) *If $n = 2t$, $|j_{n+1}(s_1; D_n)| = |j_{n+1}(s_2; D_n)| = 2^{3t}$.*

PROOF. We prove the lemma by induction. It is trivial to verify that the lemma holds for $n = 1, 2$. Assume the lemma holds for $n = k - 1$. Consider $n = k$. We have $s_1 = \{1, x_{k+1}\} \cup e$ and $s_2 = \{2, x_{k+1}\} \cup e$, where $e \subset \{x_2, \ldots, x_k\}$ with $x_i = 2i - 1$ or $2i$ for $i = 2, \ldots, k$.

First, consider $x_{k+1} = 2k + 1$. By Lemma 1(a) and (b),

(6) $\qquad j_{k+1}(s_1; D_k) = 2j_k(\{1\} \cup e; D_{k-1}) + 2j_k(\{2\} \cup e; D_{k-1})$,

(7) $\qquad -j_{k+1}(s_2; D_k) = 2j_k(\{1\} \cup e; D_{k-1}) - 2j_k(\{2\} \cup e; D_{k-1})$,

where $D_{k-1}$ is the $2^{2k-2} \times 2k$ design generated by $G_{k-1} = (\mathbf{1}, I_{k-1})$.

If $n = k = 2t+1$, the assertion of $k-1 = 2t$ implies that $|j_k(\{1\} \cup e; D_{k-1})| = |j_k(\{2\} \cup e; D_{k-1})| = 2^{3t}$. Then, from (6) and (7), we conclude that either $|j_{k+1}(s_1; D_k)|$, or $|j_{k+1}(s_2; D_k)|$ must be 0 and the other must be $2^{3t+2}$.

If $n = k = 2t+2$, the assertion of $k-1 = 2t+1$ implies that either $|j_k(\{1\} \cup e; D_{k-1})|$, or $|j_k(\{2\} \cup e; D_{k-1})|$ must be 0 and the other must be $2^{3t+2}$. Then, (6) and (7) together yield $|j_{k+1}(s_1; D_k)| = |j_{k+1}(s_2; D_k)| = 2^{3t+3}$. This proves the results for $x_{k+1} = 2k + 1$.

The proof for $x_{k+1} = 2k + 2$ is similar. Therefore, the lemma holds for $n = k$. The proof is completed by induction. □

The next result describes the complete and partial words of $D_n$ and their aliasing indexes.

LEMMA 4. *Suppose that $v_n$ consists of $p$ 1's followed by $q$ 2's, where $p + q = n$. For $l = 1, 2$, let $s_l = \{l, x_2, \ldots, x_{p+1}, 2p+3, 2p+4, \ldots, 2n+2\}$ where $x_i = 2i - 1$ or $2i$ for $i = 2, \ldots, p + 1$:*

(a) *If $p = 2t + 1$, either $\rho_k(s_1; D_n)$ or $\rho_k(s_2; D_n)$ is 0 and the other is $2^{-t}$ where $k = p + 2q + 1$;*
(b) *If $p = 2t$, $\rho_k(s_1; D_n) = \rho_k(s_2; D_n) = 2^{-t}$ where $k = p + 2q + 1$;*
(c) *$\rho_k(s_0; D_n) = 1$ where $s_0 = \{1, 2, \ldots, 2p + 2\}$ and $k = 2p + 2$;*
(d) *$\rho_k(s; D_n) = 0$ for $s$ other than $s_1$, $s_2$ or $s_0$ considered in* (a), (b) *and* (c).

PROOF. (a) and (b), when $q = 0$, it follows from Lemma 3. When $q > 0$, recursively applying Lemma 2(a) or (b) yields the result.
(c) It follows from Lemmas 1(c) and 2(c).
(d) It follows from Lemmas 1(d) and 2(d). □

Now, consider half-fractions of $D_n$. Suppose that one of the last two columns of $D_n$ is chosen as the branching column. Let $D'_n$ be the resulting $2^{2n-1} \times (2n + 1)$ design.

When the last component of $v_n$ is 1 and the last column of $D_n$ is chosen as the branching column, following (4), we can write $D'_n$ as

$$(8) \qquad D'_n = \begin{pmatrix} a & b & E & \mathbf{1} \\ -b & a & E & -\mathbf{1} \end{pmatrix}.$$

The following lemma expresses the $J$-characteristics of $D'_n$ in terms of that of $D_{n-1} = (a, b, E)$.

LEMMA 5. *Suppose that the last component of $v_n$ is 1, and the last column of $D_n$ is chosen as the branching column. For any subset $e \subset \{3, 4, \ldots, 2n\}$:*
(a) *$j(\{1\} \cup e; D'_n) = -j(\{2, 2n + 1\} \cup e; D'_n) = j(\{1\} \cup e; D_{n-1}) - j(\{2\} \cup e; D_{n-1})$;*
(b) *$j(\{2\} \cup e; D'_n) = j(\{1, 2n + 1\} \cup e; D'_n) = j(\{1\} \cup e; D_{n-1}) + j(\{2\} \cup e; D_{n-1})$;*
(c) *$j(\{1, 2, 2n + 1\} \cup e; D'_n) = 2j(\{1, 2\} \cup e; D_{n-1})$;*
(d) *$j(s \cup e; D'_n) = 0$ for $s = \{1, 2\}$, or $\{2n + 1\}$.*

It is easy to verify that choosing the second last column of $D_n$ as the branching column yields a design that is equivalent to $D'_n$ in (8).

When the last component of $v_n$ is 2 and the last (or second last) column of $D_n$ is chosen as the branching column, following (5), we can write $D'_n$ as

$$(9) \qquad D'_n = \begin{pmatrix} a & b & E & \mathbf{1} \\ -a & -b & E & -\mathbf{1} \end{pmatrix}.$$

We can also express the $J$-characteristics of $D'_n$ in terms of $D_{n-1}$.

LEMMA 6. *Suppose that the last component of $v_n$ is 2 and the last column of $D_n$ is chosen as the branching column. For any subset $e \subset \{3, 4, \ldots, 2n\}$:*

(a) $j(\{1, 2n+1\} \cup e; D'_n) = 2j(\{1\} \cup e; D_{n-1})$;
(b) $j(\{2, 2n+1\} \cup e; D'_n) = 2j(\{2\} \cup e; D_{n-1})$;
(c) $j(\{1, 2\} \cup e; D'_n) = 2j(\{1, 2\} \cup e; D_{n-1})$;
(d) $j(s \cup e; D'_n) = 0$ for $s = \{1\}, \{2\}, \{2n+1\}$, or $\{1, 2, 2n+1\}$.

4.2. *Proofs of theorems.*

PROOF OF THEOREM 1. Without loss of generality, assume that $v$ consists of $p$ 1's followed by $q$ 2's, where $p + q = n$. Lemma 4 suggests that all possible words are in forms of $s_1, s_2$ or $s_0$. If $p = 2t + 1$, by Lemma 4(a), there are $2^p$ words of length $p + 2q + 1$ with aliasing index $\rho = 2^{-t}$. If $p = 2t$, by Lemma 4(b), there are $2^{p+1}$ words of length $p + 2q + 1$ with aliasing index $\rho = 2^{-t}$. By Lemma 4(c), there is 1 complete word of length $2p + 2$. This completes the proof. $\square$

PROOF OF THEOREM 2. Without loss of generality, assume that $v$ consists of $p$ 1's and $q$ 2's, where $p + q = n$.

(a) We prove the result by induction on $p$. The result is trivial when $p = 0$. Assume that it is true for $p = k - 1$. Consider $p = k$. As in (4), we can write $D_n = D_{k+q}$, where $a$ and $b$ are the balanced two-level columns and $E$ is a full factorial with $2k + 2q - 2$ columns. We need to show that $D_{k+q}$ has projectivity $2k + 1$. Consider any subset $s$ with $2k + 1$ columns of $D_{k+q}$. There are three possible cases:

(i) Both of the last two columns of $D_{k+q}$ belong to $s$. Denote $E_1 = (a, b, E)$, $E_2 = (b, -a, E)$, $E_3 = (-a, -b, E)$ and $E_4 = (-b, a, E)$. Clearly the $E_i$'s are isomorphic to each other. The assertion of $p = k - 1$ implies that each $E_i$ has projectivity $2k - 1$. Then, the projection onto $s$ contains a full $2^{2k+1}$ factorial;

(ii) None of the last two columns of $D_{k+q}$ belong to $s$. Observe that $E$ is a full factorial with $2k + 2q - 2 \geq 2k$ columns. It is easy to verify that the projection onto $s$ contains a full $2^{2k+1}$ factorial, whether $s$ includes none, one or both of the first two columns;

(iii) One of the last two columns of $D_{k+q}$ belongs to $s$ and the other does not. Observe that the projection onto the subset consisting of the first two and the last two columns has resolution $\geq 4$ and projectivity $\geq 3$. Further, observe that $E$ is a full factorial. Then, it is easy to verify that the projection onto $s$ contains a full $2^{2k+1}$ factorial whether $s$ includes none, one or both of the first two columns.

The three cases together suggest that $D_{k+q}$ has projectivity $2k + 1$. By induction, the proof is completed.

(b) The proof is similar to (a) and omitted.

□

PROOF OF THEOREM 3. Without loss of generality, assume that $v$ consists of $p$ 1's and $q$ 2's, where $p + q = n$. Let $v_{n-1}$ be the vector consisting of the first $n-1$ components of $v$, and let $D_{n-1}$ be the binary design generated by $G_{n-1} = (v_{n-1}, I_{n-1})$.

(a) Without loss of generality, assume that the last component of $v$ is 1 and that the last column of $D$ is chosen as the branching column. If $p = 2t$, by Lemma 4(a), $D_{n-1}$ has $2^{p-1}$ words of length $p + 2q$ with aliasing index $2^{-(t-1)}$. By Lemma 5(a) and (b), these $2^{p-1}$ words in $D_{n-1}$ generate $2^p$ words of length $p + 2q$ and $2^p$ words of length $p + 2q + 1$ with aliasing index $\rho = 2^{-t}$ in $D'$. If $p = 2t + 1$, by Lemma 4(b), $D_{n-1}$ has $2^p$ words of length $p + 2q$ with aliasing index $2^{-t}$. By Lemma 5(a) and (b), these $2^p$ words in $D_{n-1}$ generate $2^{p-1}$ words of length $p + 2q$ and $2^{p-1}$ words of length $p + 2q + 1$ with aliasing index $\rho = 2^{-t}$ in $D'$. So, in both cases, $D'$ has $1/\rho^2$ words of length $p + 2q = k_1 - 1$ and $1/\rho^2$ words of length $k_1$ with aliasing index $\rho = 2^{-\lfloor p/2 \rfloor}$. By Lemma 4(c), $D_{n-1}$ has 1 complete word of length $2p$, which generates a complete word of length $2p + 1$ in $D'$ by Lemma 5(c). This completes the proof.

(b) Without loss of generality, assume that the last component of $v$ is 2 and that the last column of $D$ is chosen as the branching column. By Theorem 1, $D_{n-1}$ has 1 complete word of length $2p + 2$ and $2/\rho^2$ words of length $p + 2(q-1) + 1$ with aliasing index $\rho = 2^{-\lfloor p/2 \rfloor}$. By Lemma 6(a) and (b), each partial word in $D_{n-1}$ generates a partial word of length $p + 2q = k_1 - 1$ in $D'$ with aliasing index $\rho$. Lemma 6(c) implies that the complete word in $D_{n-1}$ produces a complete word with the same length $2p + 2 = k_2$ in $D'$. This completes the proof. □

PROOF OF THEOREM 4. Without loss of generality, assume that $v$ consists of $p$ 1's and $q$ 2's, where $p + q = n$.

(a) By Theorem 2(a), $D$ has projectivity $2p + 1$. It is obvious that any half-fraction of $D$ has projectivity $\geq 2p$. By Theorem 3(a), $D'$ has a complete word of length $2p + 1$, so its projectivity is $2p$.

(b) As in (a), by Theorem 2(b), $D$ has projectivity $2p - 1$, so $D'$ has projectivity $\geq 2p - 2$.

(c) Without loss of generality, we write $D'$ as (9), where $a$ and $b$ are balanced two-level columns and $E$ is a full factorial with $2p + 2q - 2$ columns. By Theorem 2(a), $D_{n-1} = (a, b, E)$ has projectivity $2p + 1$, so is $(-a, -b, E)$. Then, it is clear that $D'$ has projectivity $2p + 1$.

(d) By Theorem 2, $D$ has projectivity $2p + 1$, so $D'$ has projectivity $\geq 2p$. □

PROOF OF THEOREM 5. Without loss of generality, we assume $f_0 = f_3 = 0$. Then, $f_2 = n - f_1$, $k_1 = f_1 + 2f_2 + 1 = 2n - f_1 + 1$ and $k_2 = 2f_1 + 2$.

According to Theorem 1 and Corollary 2, we need to consider whether the condition $k_1 \geq k_2$ holds. It is obvious that the condition $k_1 \geq k_2$ is equivalent to $f_1 \leq (2n-1)/3$. If $k_1 \geq k_2$, the resolution is $k_2 = 2f_1 + 2$, so we shall maximize $k_2$ and choose $f_1 = \lfloor (2n-1)/3 \rfloor$, since $f_1$ is an integer. If $k_1 < k_2$, the resolution is $k_1 + 1 - \rho = 2n - f_1 + 2 - \rho$, so we shall maximize $k_1$ and choose $f_1 = \lfloor (2n+1)/3 \rfloor$, which is the smallest integer that is greater than $(2n-1)/3$.

(a) When $n = 3k - 1$, the first choice leads to $f_1 = 2k - 1$, $f_2 = k$, $k_1 = k_2 = 4k$ and $R(D) = 4k$, while the second choice leads to $f_1 = 2k$, $f_2 = k - 1$, $k_1 = 4k - 1$, $k_2 = 4k + 2$ and $R(D) = 4k - 2^{-k}$. Therefore, the first choice leads to a maximum resolution design.

(b) When $n = 3k$, the first choice leads to $f_1 = 2k - 1$, $f_2 = k + 1$, $k_1 = 4k + 2$, $k_2 = 4k$ and $R(D) = 4k$, while the second choice leads to $f_1 = 2k$, $f_2 = k$, $k_1 = 4k + 1$, $k_2 = 4k + 2$ and $R(D) = 4k + 2 - 2^{-k}$. Therefore, the second choice leads to a maximum resolution design.

(c) When $n = 3k + 1$, the first choice leads to $f_1 = 2k$, $f_2 = k + 1$, $k_1 = 4k + 3$, $k_2 = 4k + 2$ and $R(D) = 4k + 2$, while the second choice leads to $f_1 = 2k + 1$, $f_2 = k$, $k_1 = 4k + 2$, $k_2 = 4k + 4$ and $R(D) = 4k + 3 - 2^{-k}$. Therefore, the second choice leads to a maximum resolution design. □

PROOF OF THEOREM 6. Note that the minimum aberration design must maximize the integer part of the resolution. As explained in the proof of Theorem 5, we only need to consider two choices: $f_1 = \lfloor (2n-1)/3 \rfloor$ or $f_1 = \lfloor (2n+1)/3 \rfloor$.

(a) When $n = 3k - 1$, the first choice leads to a minimum aberration design with $f_1 = 2k - 1$, $f_2 = k$, $k_1 = k_2 = 4k$ and $A_{4k}(D) = 3$.

(b) When $n = 3k$, the second choice leads to a minimum aberration design with $f_1 = 2k$, $f_2 = k$, $k_1 = 4k + 1$, $k_2 = 4k + 2$, $A_{4k+1}(D) = 2$ and $A_{4k+2}(D) = 1$.

(c) When $n = 3k + 1$, the first choice leads to $f_1 = 2k$, $f_2 = k + 1$, $k_1 = 4k + 3$, $k_2 = 4k + 2$, $A_{4k+2}(D) = 1$ and $A_{4k+3}(D) = 2$, while the second choice leads to $f_1 = 2k + 1$, $f_2 = k$, $k_1 = 4k + 2$, $k_2 = 4k + 4$, $A_{4k+2}(D) = 2$, and $A_{4k+4}(D) = 1$. Therefore, the first choice leads to a minimum aberration design. □

PROOF OF THEOREM 7. It follows from Theorem 2. □

PROOF OF THEOREM 8. Without loss of generality, we assume $f_0 = f_3 = 0$ so that $f_1 + f_2 = n$. According to Theorem 3, we need to consider four cases: (i) the branching column is associated with number 1 and $k_1 \geq k_2$, (ii) the branching column is associated with number 1 and $k_1 < k_2$, (iii) the branching column is associated with number 2 and $k_1 - 1 \geq k_2$ and (iv) the branching column is associated with number 2 and $k_1 - 1 < k_2$. For each case, we choose $f_1$ and $f_2$ to maximize the shortest wordlength and resolution. The resolutions and wordlength patterns of the resulting designs can be calculated by Corollaries 4, 5, 6 and 7.

(a) When $n = 3k - 1$, the condition $k_1 \geq k_2$ is equivalent to $f_1 \leq 2k - 1$; the condition $k_1 - 1 \geq k_2$ is equivalent to $f_1 \leq 2k - 4/3$. For case (i), we want to maximize $k_2$, so we choose $f_1 = 2k - 1$ and $f_2 = k$, which yields $k_1 = 4k$, $k_2 = 4k$, $R(D') = 4k - 1$, $A_{4k-1}(D') = 2$ and $A_{4k}(D') = 1$. For case (ii), we want to maximize $k_1$, so we choose $f_1 = 2k$ and $f_2 = k - 1$, which yields $k_1 = 4k - 1$, $k_2 = 4k + 2$, $R(D') = 4k - 1 - 2^{-k}$, and $A_{4k-2}(D') = A_{4k-1}(D') = A_{4k+1}(D') = 1$. For case (iii), we want to maximize $k_2$, so we choose $f_1 = 2k - 2$ and $f_2 = k + 1$, which yields $k_1 = 4k + 1$, $k_2 = 4k - 2$, $R(D') = 4k - 2$, $A_{4k-2}(D') = 1$ and $A_{4k}(D') = 2$. For case (iv), we want to maximize $k_1$, so we choose $f_1 = 2k - 1$ and $f_2 = k$, which yields $k_1 = 4k$, $k_2 = 4k$, $R(D') = 4k - 2^{-(k-1)}$, $A_{4k-1}(D') = 2$ and $A_{4k}(D') = 1$. Therefore, the design in case (iv) has both maximum resolution and minimum aberration.

(b) When $n = 3k$, the condition $k_1 \geq k_2$ is equivalent to $f_1 \leq 2k - 1/3$; the condition $k_1 - 1 \geq k_2$ is equivalent to $f_1 \leq 2k - 2/3$. For case (i), we shall choose $f_1 = 2k - 1$ and $f_2 = k + 1$, which yields $k_1 = 4k + 2$, $k_2 = 4k$, $R(D') = 4k - 1$ and $A_{4k-1}(D') = A_{4k+1}(D') = A_{4k+2}(D') = 1$. For case (ii), we shall choose $f_1 = 2k$ and $f_2 = k$, which yields $k_1 = 4k + 1$, $k_2 = 4k + 2$, $R(D') = 4k + 1 - 2^{-k}$, $A_{4k}(D') = 1$ and $A_{4k+1}(D') = 2$. For case (iii), we shall choose $f_1 = 2k - 1$ and $f_2 = k + 1$, which yields $k_1 = 4k + 2$, $k_2 = 4k$, $R(D') = 4k$, $A_{4k}(D') = 1$ and $A_{4k+1}(D') = 2$. For case (iv), we shall choose $f_1 = 2k$ and $f_2 = k$, which yields $k_1 = 4k + 1$, $k_2 = 4k + 2$, $R(D') = 4k + 1 - 2^{-k}$, $A_{4k}(D') = 2$ and $A_{4k+2}(D') = 1$. Therefore, the design in case (ii) has both maximum resolution and minimum aberration.

(c) When $n = 3k + 1$, the condition $k_1 \geq k_2$ is equivalent to $f_1 \leq 2k + 1/3$; the condition $k_1 - 1 \geq k_2$ is equivalent to $f_1 \leq 2k$. For case (i), we shall choose $f_1 = 2k$ and $f_2 = k + 1$, which yields $k_1 = 4k + 3$, $k_2 = 4k + 2$, $R(D') = 4k + 1$ and $A_{4k+1}(D') = A_{4k+2}(D') = A_{4k+3}(D') = 1$. For case (ii), we shall choose $f_1 = 2k + 1$ and $f_2 = k$, which yields $k_1 = 4k + 2$, $k_2 = 4k + 4$, $R(D') = 4k + 2 - 2^{-k}$ and $A_{4k+1}(D') = A_{4k+2}(D') = A_{4k+3}(D') = 1$. For case (iii), we shall choose $f_1 = 2k$ and $f_2 = k + 1$, which yields $k_1 = 4k + 3$, $k_2 = 4k + 2$, $R(D') = 4k + 2$ and $A_{4k+2}(D') = 3$. For case (iv), we shall choose $f_1 = 2k + 1$ and $f_2 = k$, which yields $k_1 = 4k + 2$, $k_2 = 4k + 4$, $R(D') = 4k + 2 - 2^{-k}$, $A_{4k+1}(D') = 2$ and $A_{4k+4}(D') = 1$. Therefore, the design in case (iii) has both maximum resolution and minimum aberration. □

PROOF OF THEOREM 9. It follows from Theorem 4. □

PROOF OF THEOREM 10. (a) The quarter-fraction designs given in Theorems 6 and 8 have the same wordlength patterns as the regular minimum aberration designs. Then, the result follows from Theorem 2 of Xu [27], which states that the regular minimum aberration $2^{m-2}$ design has minimum aberration among all possible designs.

(b) The quarter-fraction designs given in Theorems 7 and 9 have $2^{m-2}$ runs and projectivity $m - 3$. It is sufficient to prove that the projectivity of any $2^k \times m$ two-level design $D$ is at most $k - 1$ when $m \geq k + 2$. Assume that $D$ has projectivity $k$. Then, the projection onto any $k$ factors is an unreplicated $2^k$ full factorial, because $D$ has exactly $2^k$ runs. Therefore, $D$ is an orthogonal array of strength $k$. Theorem 2.19 of Hedayat, Sloane and Stufken [15] implies that $m < k + 1$. This contradicts the condition $m \geq k + 2$. □

Department of Statistics
University of California
Los Angeles, California 90095–1554
USA
E-mail: fredphoa@stat.ucla.edu
hqxu@stat.ucla.edu